\documentclass[12pt, oneside]{article}
\usepackage{amsmath,amsfonts,amssymb,amscd,theorem}
\date{}

\newtheorem{thm}{Theorem}

\newtheorem{lem}{Lemma}

\begin{document}

 \title{K\"ahler manifolds and fundamental groups
of negatively $\delta$-pinched manifolds }

\author{J\"urgen Jost and Yi-Hu Yang\thanks{The second named author
supported partially by NSF of China (No. 10171077)}}

 \maketitle

\begin{abstract}
In this note, we will show that the fundamental group of any
negatively $\delta$-pinched ($\delta > {\frac 1 4}$) manifold
can't be the fundamental group of a quasi-compact K\"ahler
manifold. As a consequence of our proof, we also show that any
nonuniform lattice in $F_{4(-20)}$ cannot be the fundamental group
of a quasi-compact K\"ahler manifold. The corresponding result for
uniform lattices was proved by Carlson and Hern\'andez \cite{ch}.
Finally, we follow Gromov and Thurston \cite{gt} to give some
examples of negatively $\delta$-pinched manifolds ($\delta >
{\frac 1 4}$) of finite volume which, as topological manifolds,
admit no hyperbolic metric with finite volume under any smooth
structure. This shows that our result for $\delta$-pinched
manifolds is a nontrivial generalization of the fact that no
nonuniform lattice in $SO(n, 1) (n \ge 3)$ is the fundamental
group of a quasi-compact K\"ahler manifold \cite{y}.
 \end{abstract}

\section{Introduction}

In \cite{yz}, Yau and Zheng (independently, Hern\'andez \cite{h})
studied negatively $\delta (\ge {\frac 1 4})$-pinched manifolds.
In particular, they showed that such manifolds are hermitian
negative (for the definition, see Section 2). An interesting
consequence of their result (which was not stated explicitly in
\cite{yz}, but implied clearly) is that the fundamental group of a
compact negatively $\delta (>{\frac 1 4})$-pinched manifold cannot
be the fundamental group of any compact K\"ahler manifold (in the
sequel, we always assume $\delta > {\frac 1 4}$). On the other
hand, in \cite{gt}, Gromov and Thurston constructed some examples
of closed $n (n\ge 4)$-manifolds which admit some negatively
$\delta$-pinched metrics, but no metric of constant curvature $-1$
under any smooth structure. In other words, $\pi_1$ of these
manifolds can't be a cocompact discrete subgroup of $SO(n, 1)$. A
simple homotopical (or cohomological dimension) argument also
shows that these $\pi_1$ cannot be cocompact discrete subgroups of
other $SO(m, 1), m\neq n$. So, the above result by Yau-Zheng and
Hern\'andez is a nontrivial generalization of a result by
Carlson-Toledo \cite{ct} and Jost-Yau \cite{jy1} independently (in
\cite{jy1}, although the result was not stated explicitly, it is
clear that it is contained in the results of \cite{jy1}), which
asserts that no cocompact lattice in $SO(n,1) (n\ge 3)$ can be
$\pi_1$ of a compact K\"ahler manifold. In \cite{y}, the second
author also showed that the result by Carlson-Toledo and Jost-Yau
is still valid for nonuniform lattices: Let $\overline{M}$ be a
compact K\"ahler manifold and $D$ be a normal crossing divisor,
denote ${\overline M}\setminus D$ by $M$. Call $M$ a quasi-compact
K\"ahler manifold. Topologically, this class of manifolds includes
quasi-projective varieties by Hironaka's Theorem for resolution of
singularities. Then, no nonuniform lattice in $SO(n, 1) (n\ge 3)$
can be $\pi_1$ of a quasi-compact K\"ahler manifold. In this note,
our purpose is to generalize this to open $-\delta$-pinched
manifolds with finite volume, namely we will show the following

\begin{thm}
Let $M$ be a quasi-compact K\"ahler manifold and $N$ a complete
noncompact $-\delta$-pinched manifold with finite volume and
dimension $\ge 3$. Then $\pi_1(N)$ is not isomorphic to
$\pi_1(M)$. Namely, the fundamental group of a complete noncompact
$-\delta$-pinched manifold with finite volume and dimension $\ge
3$ cannot be the fundamental group of any quasi-compact K\"ahler
manifold.
 \end{thm}

The method of proof is to use harmonic map theory due to Jost and
Zuo \cite{jz0, jz}. One of the ingredients for the proof is to
show that a differentiable family of submanifolds of codimension
$2$ from the harmonic map in question actually is a holomorphic
family and it gives rise to a foliation. Also observe that the
argument of Theorem 7.1 in \cite{ct} does not work in the present
noncompact situation. In addition, we need to treat the case of
dimension $3$ separately. As an interesting corollary of the above
theorem, we have
\begin{thm}
Let $M_1, M_2, \cdots , M_s (s\ge 2)$ be some compact or open
Riemann surfaces. Then, the product $M_1\times
M_2\times\cdots\times M_s$ does not admit any complete negatively
$\delta$ ($\delta > {\frac 1 4}$)-pinched metric with finite
volume under any smooth structure.
\end{thm}

This is analogous to the positive curvature case, where by the
$\frac 1 4$-sphere theorem, one knows that $S^2\times
S^2\times\cdots\times S^2$ (at least two copies) does not admit
any $\delta$-pinched metric.

Combining the technique of the proof of Theorem 1 with the
Lie-theoretic analysis of \cite{ch}, one can also obtain
\begin{thm}
Let $\Gamma$ be a nonuniform lattice in the exceptional group
$F_{4(-20)}$. Then $\Gamma$ can't be the fundamental group of any
quasi-compact K\"ahler manifold.
\end{thm}

\vskip 4mm \noindent {\bf Remark:} For the cocompact lattice case
of $F_{4(-20)}$, the above result was proved by Carlson and
Hern\'andez \cite{ch}.

\vskip 4mm \noindent {\bf Acknowledgements:} The second named
author would like to thank Kang Zuo for his constant help. This
research began when the second named author was visiting the
Max-Planck-Institute for Mathematics in the Sciences in Leipzig.
He thanks the institute for its hospitality and good working
conditions. Finally, the authors would like to thank a
referee for  pointing out that  initial examples
were incorrect in the first version of the paper.

\section{Negatively $\delta$-pinched manifolds and harmonic maps}

In this section, we will recall some well-known results concerning
negatively $\delta$-pinched manifolds \cite{bk, yz, h} and the
existence of harmonic maps \cite{jz} and give some properties of
the harmonic maps in question.

Let $N$ be a complete negatively $\delta$-pinched manifold with
finite volume. By means of the Margulis lemma, one gets that as
stated and proved in \cite{bk}, Corollary 1.5.2, $N$ can be
considered as the interior of a compact manifold with boundary,
and the boundary topologically is the disjoint union of tori up to
a finite group action. And concerning the curvature tensor of $N$,
one has

\begin{lem} (\cite{yz}, Lemma 2, 3 or \cite{h}, Theorem 2.5)
$N$ is Hermitian negative, namely, for any two complex vectors $Z,
W$ in $T_p(N)\otimes {\mathbb{C}}$ ($p\in N$),
\[
R(Z, W, \overline{Z}, \overline{W}) \le 0;
\]
and there do not exist complex linearly independent vectors $Z, W$
satisfying
\[
R(Z, W, \overline{Z}, \overline{W}) = 0.
\]
\end{lem}

Now, we turn to harmonic map theory. Let $\overline M$ be a
compact K\"ahler manifold of dimension $n$ with a fixed K\"ahler
metric $\omega_0$, $D$ be a fixed divisor with (at worst) normal
crossings and $D=\bigcup_{i=1}^{p}D_{i}$. Here, $D_i$ are the
irreducible components of $D$. One may also assume that each
irreducible component $D_i$ is free from self intersections. Thus,
at each intersection point, precisely at most $n$ components of
$D$ meet. Denote $\overline M\setminus D$ by $M$.

Let $\sigma_i ~(i=1, 2, \cdots , p)$ be a defining section of
$D_i$ in $\mathcal{O}({\overline M}, [D_i])$, which satisfies
$|\sigma_i|\le 1$ for a certain Hermitian metric of $[D_i]$ and
defines a holomorphic coordinate system in each small disk
transversal to $D_i$. So, one can get a fibration of a small
neighborhood, say $|\sigma_i|\le\mu\le 1$, of $D_i$ by small
holomorphic disks which meet $D_i$ transversally. Similarly, for
the boundary of such a small neighborhood, denoted by
$\Sigma_i^{\mu}$, one also gets a fibration by circles. The
intersection of two such boundaries is fibered by tori.

Corresponding to the above defining sections $\sigma_i$, one can
define a complete K\"ahler metric on $M$ as follows,
\[
   g:=-{\frac{\sqrt{-1}}2}\sum_{i=1}^p\partial\overline{\partial}(\phi
(|\sigma_i|){\text{log}}|{\text{log}}|\sigma_i|^2|)+c\omega_0|_M,
\]
where $\phi$ is a suitable $C^{\infty}$ cut-off function on $[0,
\infty )$, so that $\phi (s)$ is identical to one on $[0, \epsilon
)$ and to zero on $[2\epsilon , \infty )$, for sufficiently small
$\epsilon\ge 0$, and $c$ is taken sufficiently large, so that $g$
is positive definite. Then $g$ is a K\"ahler metric. One can show
that $(M, g)$ is complete and has finite volume \cite{cg}. In
fact, when restricted to small holomorphic disks transversal to
$D$, this metric looks like the Poincar\'e metric on the punctured
disk $(D^*, z)$
\[
-{\frac{\sqrt{-1}}{2}}\partial\overline{\partial}{\text{log}}
(-{\text{log}}|z|^2)={\frac{\sqrt{-1}}{2}}{\frac{{\text d}z\wedge
{\text d}{\overline z}}{|z|^2({\text{log}}|z|^2)^2}} .
\]
In the sequel, unless stated otherwise, we always assume that $M$
is endowed with the above complete metric $g$.

As before, let $N$ be a complete noncompact negatively
$\delta$-pinched manifold with finite volume, whose universal
covering is denoted by $\widetilde{N}$. Denote the isometry group
of $\widetilde{N}$ by $I(\widetilde{N})$. Given a {\it reductive}
homomorphism (for the definition of reductivity, see the
Definition 1.1 of \cite{jz})
\[
\rho : \pi_1(M)\to I(\widetilde{N}),
\]
one wants to get a $\rho$-equivariant harmonic map from the
universal covering of $M$ to $\widetilde{N}$. In general,
difficulties will arise since the homomorphism $\rho$ may map some
small loops around $D$ to some {\it hyperbolic} or {\it
quasi-hyperbolic} elements (for their definitions and the
definitions of {\it elliptic} and {\it parabolic elements}, see
the Definition 1.2 of \cite{jz}) in $I(\widetilde{N})$. This is
why a $\rho$-equivariant harmonic map, if it exists, may have
infinite energy (here, we use the above metric $g$ to compute the
energy). It is however  fortunate that this 
difficulty will not arise for the application in this paper. In
the following, we will assume that the image of $\rho$ lies in
$\pi_1(N)$. (Actually, in the application, $\rho$ is an
isomorphism from $\pi_1(M)$ to $\pi_1(N)$.) 

As in \cite{jz}, one needs to consider only two cases: 1) $\rho$
maps every small loop around a (topological) component of $D$ (as
an element in $\pi_1(M)$) to a hyperbolic or elliptic element; 2)
$\rho$ maps every small loop around a component of $D$ to a
parabolic or quasi-hyperbolic element. It is useful to point out
that any two loops in each component of $D$, as elements of the
fundamental group, commute with each other, so $\rho$ maps them
simultaneously to either hyperbolic (elliptic) elements or
quasihyperbolic (parabolic) elements. Since $N$ is a complete
noncompact negatively $\delta$-pinched manifold with finite
volume, so the image elements can not be quasihyperbolic
(otherwise, $N$ will be of infinite volume) and therefore only
parabolic images may occur in the second case. For the parabolic
images case, the problem can be handled as done in pages 85-91 of
\cite{jz0}, where all steps can be translated to the present
situation after knowing the fact that the Jacobi fields in the
present case are also exponentially decaying, which is an easy
consequence of Riemannian geometry. We now  address 
the first case. Similar to the parabolic images case, the elliptic
images case will  also not cause any difficulty and it can also be
handled as in \cite{jz0}. Thus, one only needs to handle the case
of hyperbolic images. Since $N$ is a negatively $\delta$-pinched
manifold ($\delta >{\frac 1 4}$), so the situation is similar to
that of locally symmetric spaces of noncompact type of rank one
and slightly simpler than in the situation of general symmetric
spaces. In this case, one furthermore has the fact that the
element of $\pi_1(\Sigma^{\mu})$ (here $\Sigma^{\mu}$ is the boundary
of the $\mu$-tube neighborhood of $D$) determined by any fibre of
$\Sigma^{\mu}_{i}\to D_{i}$ is central and hence independent of the
choice of the fibre and the fact that the centralizer of a
hyperbolic element in a discrete subgroup of $I(\widetilde{N})$ is
virtually cyclic. Using these facts, one can then follow the
arguments in pages 477-480 of \cite{jz}  to get a
desired $\rho$-equivariant map and hence a $\rho$-equivariant
harmonic map. Finally, the argument in pages 481-483 of \cite{jz}
(Lemma 1.1) shows that the harmonic map is pluriharmonic.

Now, we can state the following theorem on the existence of a
$\rho$-equivariant pluriharmonic map, which is essentially due to
Jost and Zuo \cite{jz0, jz}, as follows:

\begin{thm}
Let $M$, $N$, $I(\widetilde{N})$ and $\rho$ as above. Then there
exists a $\rho$-equivariant pluriharmonic map $u$ from the
universal covering $\widetilde{M}$ of $M$ with the above metric
$g$ to $\widetilde{N}$.
\end{thm}

\noindent {\bf Remark}: Actually, Jost-Zuo's theorem in \cite{jz0}
and \cite{jz} is more general in some respects. In particular, if
$N$ is a locally symmetric space of noncompact type of rank one
with finite volume, the theorem still holds, as will be used in
the proof of the Theorem 3.

For simplicity of notation, we shall consider the harmonic map $u$
in the sequel as a map from $M$ to $N' = \widetilde{N}/\rho
(\pi_1(M))$. (In our applications $N'$ is just $N$). Let $w\in
D_i$ be a regular point of $D$. Near $w$, one can choose a
coordinate system $(z^1, z^2)$ on $\overline M$ such that $z^1$
parameterizes small holomorphic discs, which meet $D_i$
transversally near $w$, $z^2$ parameterizes $D_i$ (of course,
$z^2$ will have more than one component if the complex dimension
of $M$ is greater than $2$. In the following, the index $2$ will
stand for all those $z^2$-directions together), and $z^1=0$ on a
small neighborhood of $w$ in $D_i$ and $z^2(w)=0$. One then has
some derivative estimates for $u$ (see p.481 of \cite{jz}):
\[
|{\frac{\partial u}{\partial z^1}}(z^1, z^2)|_g\le{\frac{c}{|z^1|}},~~~~
|{\frac{\partial u}{\partial z^2}}(z^1, z^2)|_g\le c,
\]
where $c$ is some positive constant. If $w$ is a singular point of
$D$, i.e.,a point at which two irreducible components of $D$ meet,
similar estimates can be obtained. One may use
$\sigma:=\prod_{i=1}^{p}\sigma_i$ to replace the above coordinate
component $z^1$. Then, one can get that in the $\sigma$-direction,
the derivative of $u$ behaves like ${\frac{1}{|\sigma |}}$,
whereas in directions normal to $\sigma$, it is bounded.

Now, we shall show that the rank of the harmonic map $u$
has a serious restriction.
Let $M$ be as before with the constructed K\"ahler metric $g$,
the corresponding K\"ahler form of which is denoted by
\[
\omega =
{\frac{\sqrt{-1}}2}\sum_{\alpha , \beta =1}^m
g_{\alpha\overline{\beta}}dz^{\alpha}\wedge dz^{\overline{\beta}}
\]
where $m=\text{dim}_{\mathbb{C}}M$ and $(z^1, z^2, \cdots , z^m)$
is a local coordinate system of $M$. Introduce a local coordinate
system $(u^1, u^2, \cdots , u^n)$ on $N'$. As in \cite{s1}, we
introduce a symmetric $(2, 0)$-tensor $\phi$ related to the map
$u$
\[
\phi (X, Y)=<\partial u(X), \partial u(Y)>,~~~~X, Y\in T^{1,0}_xM,
\]
which can be locally written as $\sum_{\alpha,\beta
=1}^{m}\phi_{\alpha\beta}{\text d}z^{\alpha}\otimes{\text
d}z^{\beta}$. Now, we compute its iterated divergence. By the
divergence formula, one has a $(1, 0)$-form $\xi$
\[
\xi_{\alpha}=g^{\beta\overline{\gamma}}\phi_{\alpha\beta ,\overline{\gamma}}
\]
where $(g^{\beta\overline{\gamma}})$ represents the inverse of
$(g_{\beta\overline{\gamma}})$ and "," denotes the covariant
derivative. Then, taking the divergence of $\xi$ again, one
obtains, by a direct computation \cite{s1},
\[
\delta\xi
=(|D''\partial u|^2
-g^{\alpha\overline{\beta}}g^{\gamma\overline{\delta}}
R_{iklm}u^i_{\alpha}u^k_{\gamma}u^l_{\overline{\beta}}u^m_{\overline{\delta}})
\]
where $\delta$ is the codifferential, $R_{iklm}$ is the curvature
tensor of $N'$, and $D''\partial u$ is the (0, 1)-type covariant
derivative of $\partial u$, which is locally written as
\[
(D''\partial u)^i_{\alpha\overline{\beta}}
=u^i_{\alpha\overline{\beta}}+\Gamma^i_{jk}u^j_{\alpha}u^k_{\overline{\beta}}.
\]
Here $\Gamma^i_{jk}$ are the Christoffel symbols of $Y$. Then,
Jost-Zuo's argument (see p.482 of \cite{jz}) shows that the two
sides of the above formula are zero pointwise. It should be
pointed out that the estimates given above of the derivatives of
$u$ near the divisor $D$ are essential in this reasoning. In
particular, using the curvature conditions of $N'$, one obtains
that
\[
D''\partial u=0,~~~~
g^{\alpha\overline{\beta}}g^{\gamma\overline{\delta}}
R_{iklm}u^i_{\alpha}u^k_{\gamma}u^l_{\overline{\beta}}u^m_{\overline{\delta}}=0.
\]
Note that the above first formula just represents the pluriharmonicity of $u$.
Taking the holomorphic orthogonal frame
 $e_1, e_2, \cdots , e_m $
on $M$, one has
\[
g^{\alpha\overline{\beta}}g^{\gamma\overline{\delta}}
R_{iklm}u^i_{\alpha}u^k_{\gamma}u^l_{\overline{\beta}}u^m_{\overline{\delta}} =
<R(\partial u(e_{\alpha}), \partial u(e_{\beta}))\overline{\partial}
u(e_{\overline{\alpha}}), \overline{\partial}u(e_{\overline{\beta}})>.
\]
Thus, by means of the previous lemma, there exist two complex
constants $a, b$ (at least one $\neq 0$) satisfying $a\partial
u(e_{\alpha}) + b\partial u(e_{\beta}) = 0$, namely,
$du(T^{1,0}M)$ is complex one-dimensional. So, one obtains

\begin{lem}
Let $u:M\to N'$ be the pluriharmonic map as in the Theorem $4$.
Then $u$ has real rank at most $2$.
\end{lem}

As a consequence of the above arguments, one has

\begin{thm}
Let $M, N$ be as before, and let $u: M\to N$ be a smooth map. Then
$u$ is harmonic if and only if $D''$ is the
$\overline{\partial}$-operator of a holomorphic structure on
$u^*T^{\mathbb{C}}N$ and $\partial u$ is a holomorphic section of
the bundle ${\text{Hom}}(T^{1,0}M,
\\u^*T^{\mathbb{C}}N)$.
\end{thm}

\noindent {\bf Proof:} It is easy to see that $D''$ is the
$\overline{\partial}$-operator of a holomorphic structure on
$u^*T^{\mathbb{C}}N$ if and only if it satisfies the integrability
condition $(D'')^2 = 0$. An easy computation shows
\[
(D'')^2(X, Y) = R(du(X), du(Y)),
\]
where $R$ is the complex-multilinear extension of the curvature
tensor of $N$ and $X, Y\in T^{0, 1}M$. Since $du(T^{1,0}M)$ is
complex one-dimensional, $(D'')^2(X, Y) = 0$.

\section{Factorization of harmonic maps and proofs of theorems}

From the argument of the previous section, we know that if $f$ is
a harmonic map from $M$ (a quasi-compact K\"ahler manifold with an
appropriate complete K\"ahler metric constructed as in the
previous section) to $N$ (a complete noncompact negatively
$\delta$-pinched manifold with finite volume ) from the Theorem
$4$, then it has real rank at most 2. From now on, we assume that
$\pi_1 (N)$ is isomorphic to $\pi_1 (M)$. We want to derive a
contradiction, hence the Theorem $1$ is proved.

Let $f$ be a harmonic map from $M$ to $N$ from the Theorem $4$
which induces an isomorphism from $\pi_1 (M)$ to $\pi_1 (N)$. From
the argument of the previous section, we know that $f$ has real
rank at most $2$. Obviously, its real rank cannot be zero; if it
has real rank one, by a result of J. H. Sampson \cite{s2}, $f$
maps $M$ to a closed geodesic in $N$. So, $\pi_1 (N)$ is
isomorphic to the ring of integers $\mathbb{Z}$. Because $N$ is
noncompact but of finite volume, its fundamental group must
contain some parabolic element. So, each element in $\pi_1 (N)$ is
parabolic. This is impossible, since $N$ has finite volume. So, we
can assume that $f$ has real rank $2$ at some point. Actually, it
has real rank $2$ generically, since it is pluriharmonic (more
precisely, $f^*TN\otimes {\mathbb{C}}$ is a holomorphic bundle
under the $(0,1)$-part of the induced connection and $\partial f$
is a holomorphic section of the bundle ${\text{Hom}}(T^{1,0}M,
f^*TN\otimes {\mathbb{C}})$).

In the following, we will show that $f$ gives rise to a foliation
on $M$. Note that the argument of Theorem 7.1 in \cite{ct} does
not work in the noncompact case. We will make use of a similar
(but more general) argument of \cite{jy1}. Assume $f$ has real
rank $2$ at the point $z_0$, so $f$ has real rank $2$ in a
neighborhood, say, $U$. Take a holomorphic coordinate system
$(z^1, z^2, \cdots , z^m )$ around $z_0$ (w.l.o.g., one can assume
the system covers $U$). By the previous section's result,
$df(T^{1,0}M)$ is complex one-dimensional. Without loss of
generality, one can assume $df({\frac{\partial}{\partial z^1}})
\neq 0$ and denote it by $X$. Set $\partial
f({\frac{\partial}{\partial z^i}}) = q^iX, i = 2, 3, \cdots , m$,
where $\{q^i\}$ are complex valued functions defined on $U$. We
will show that they are actually holomorphic. Since $f$ has real
rank $2$, $\bar X \neq X$ and $\bar X \neq -X$. So, for $i = 2, 3,
\cdots , m, j = 1, 2, \cdots , n$, one has, using the
K\"ahlerianity of $M$ and the pluriharmonicity of $f$
\begin{eqnarray*}
0 &=&
D''df({\frac{\partial}{\partial z^i}},
{\frac{\partial}{\partial{\bar z}^j}})\\
&=&
D''_{{\frac{\partial}{\partial{\bar z}^j}}}
df({\frac{\partial}{\partial z^i}})\\
&=&
D''_{{\frac{\partial}{\partial{\bar z}^j}}}
(df({\frac{\partial}{\partial z^i}}))\\
&=&
D''_{\frac{\partial}{\partial{\bar z}^j}}(q^iX)\\
&=&
{\frac{\partial q^i}{\partial{\bar z}^j}}X +
q^iD''_{\frac{\partial}{\partial{\bar z}^j}}
(df(\frac{\partial}{\partial z^1}))\\
&=&
{\frac{\partial q^i}{\partial{\bar z}^j}}X.
\end{eqnarray*}
So, $q^i, i = 2, 3, \cdots , m$ are holomorphic. Consider the
holomorphic distribution on $U$
\[
\{  {\frac{\partial}{\partial z^2}}-q^2{\frac{\partial}{\partial z^1}},
{\frac{\partial}{\partial z^3}}-q^3{\frac{\partial}{\partial z^1}},
\cdots ,
{\frac{\partial}{\partial z^m}}-q^m{\frac{\partial}{\partial z^1}}
\}.
\]
Obviously, it is the holomorphic kernel of the differential $df$.
Moreover, by the complex version of a standard fact (See,
\cite{kl}, Proposition 1.4.10), $df([{\frac{\partial}{\partial
z^i}} -q^i{\frac{\partial}{\partial z^1}},
{\frac{\partial}{\partial z^j}}- q^j{\frac{\partial}{\partial
z^1}}]) = 0, i, j = 2, 3, \cdots , m, i.e.,
[{\frac{\partial}{\partial z^i}} -q^i{\frac{\partial}{\partial
z^1}}, {\frac{\partial}{\partial z^j}}-
q^j{\frac{\partial}{\partial z^1}}])$ lie in this distribution.
So, it is a holomorphic integrable distribution. Then, the complex
version of Frobenius theorem asserts that on a neighborhood of
$z_0$ (assume it still is $U$), there is a holomorphic coordinate
system $(w^1, w^2, \cdots , w^m)$ satisfying
\[
\{
{\frac{\partial}{\partial w^2}}, {\frac{\partial}{\partial w^3}},
\cdots , {\frac{\partial}{\partial w^m}}
\} =
\{  {\frac{\partial}{\partial z^2}}-q^2{\frac{\partial}{\partial z^1}},
{\frac{\partial}{\partial z^3}}-q^3{\frac{\partial}{\partial z^1}},
\cdots ,
{\frac{\partial}{\partial z^m}}-q^m{\frac{\partial}{\partial z^1}}
\}.
\]
That is to say, $df({\frac{\partial}{\partial w^i}}) = 0, i=2, 3,
\cdots , m$, but $df({\frac{\partial}{\partial w^1}}) \neq 0$. So,
when restricted to the hypersurfaces $w^1 = const$, $f$ is
constant. Therefore, one obtains a well-defined foliation
$\cal{F}$ on a Zariski open set of $M$, namely, on the set where
$f$ has real rank $2$. Arguments of Mok (see Proposition (2.2.1)
of \cite{m1}) imply that $\cal{F}$ can be extended as a
holomorphic foliation to $M\setminus V$ for some complex analytic
variety $V$ of complex codimension at least $2$. Then, the study
of \cite{m2} (see Proposition (2.2) of \cite{m2}) shows that the
extended foliation actually defines an open analytic equivalence
relation, still denoted by $\cal F$, on $M$, and the quotient of
$M$ by $\cal{F}$, denoted by $S$, is an irreducible complex space
of complex dimension $1$, by a result of Kaup \cite{ka}.
Therefore, one has a factorization of the harmonic map $f$: $f =
h\circ\pi$, where $\pi : M\to S$ is holomorphic because of the
construction of $S$ and $h: S\to N$ is harmonic since $f$ is
pluriharmonic.

Since $f_{*}: \pi_1(M)\to \pi_1(N)$ is an isomorphism, $\pi_{*}:
\pi_1(M)\to\pi_1(S)$ is injective. Therefore, $\pi_1(N)$, as a
subgroup of $\pi_1(S)$, acts freely on the universal covering of
$S$, which is contractible as a topological space. Thus, the
cohomological dimension of $\pi_1(N)$ \cite{b} is at most $2$.
But, since $N$ is a negatively $\delta$-pinched manifold with
finite volume, it can topologically be regarded as the interior of
a manifold with boundary, here the boundary is the disjoint union
of tori up to a finite group. So, by means of a result of \cite{b}
(p. 211, Corollary 8.3), $\pi_1(N)$ has cohomological dimension
$n-1$, here $n$ is the real dimension of $N$. So, if $n\ge 4$, we
derive a contradiction. In the following, we will treat the case
$n=3$ separately. The idea of the proof was told to us by
Professor M. S. Raghunathan.

We now assume that $N$ has dimension $3$ and $\pi_1(N)$ is
isomorphic to $\pi_1(M)$. So, by means of the previous argument,
one has a holomorphic map $h$ from $M$ to an irreducible complex
space $S$ of complex dimension $1$, which induces an injective map
$h_{*}: \pi_1(M)\to\pi_1(S)$. Now, we have two cases to discuss:
1. $S$ is noncompact; 2. $S$ is compact.

\vskip 4mm \noindent {\bf Case 1}: If $S$ is noncompact, a
standard argument shows that $\pi_1(S)$ is free, so its
cohomological dimension is $1$, consequently, the cohomological
dimension of its subgroup is also $1$. But the cohomological
dimension of $\pi_1(N)$ is $2$. So only the case 2 may occur.

\vskip 4mm \noindent {\bf Case 2}: In this case, we also have two
cases to discuss: i) the image of $h_*$ has infinite index in
$\pi_1(S)$; ii) the image of $h_*$ has finite index in $\pi_1(S)$.
If the case i) is true, one can lift $h$ so that the case can be
actually reduced to the case 1). So, that case is impossible; if
the case ii) is true, by a finite lifting, one can also assume
that $h_{*}$ is an isomorphism from $\pi_1(M)$ to $\pi_1(S)$. We
will also derive a contradiction. Since $\pi_1(M)$ is isomorphic
to $\pi_1(N)$ by the assumption, so $\pi_1(N)$ is isomorphic to
$\pi_1(S)$. Since $N$ is of negatively $\delta$-pinched curvature
and finite volume, by means of Corollary 1.5.2 in \cite{bk}, one
can consider $N$ as the interior of a compact manifold with
boundary, here the boundary is the disjoint union of tori up to a
finite group action. Then an easy exercise shows that the
Euler-Poincar\'e characteristic $\chi (\pi_1(N))$ of $\pi_1(N)$ is
zero by using the homological exact sequence of a space pair $(X,
\partial X)$ for a compact manifold $X$ and its boundary $\partial
X$ and Poincar\'e duality theorem for manifolds with boundary (see
\cite{ma}, p. 227).

We now show that the Euler-Poincar\'e characteristic of $\pi_1(S)$
is not zero and hence a contradiction. Since $S$ is a compact
irreducible complex space of complex dimension $1$, so, by passing
to normalizations, without loss of generality we can assume that
$S$ is normal and hence smooth and the above factorization $f =
h\circ\pi$ remains valid. Clearly, $S$ can't be the sphere since
$\pi_1(N)$ acts freely on the universal covering of $S$. If $S$ is
of genus $g\ge 2$, then the Euler-Poincar\'e characteristic of
$\pi_1(S)$ is $2-2g < 0$. This is a contradiction. So, if the
assumption is true, $S$ must be a torus. Thus $\pi_1(N) =
{\mathbb{Z}} + {\mathbb{Z}}$. In the following, using harmonic map
theory, we will show that $\pi_1(N)$ cannot be ${\mathbb{Z}} +
{\mathbb{Z}}$. Assume $\pi_1(N) = {\mathbb{Z}} + {\mathbb{Z}}$.
Take the standard torus $T$ with a flat metric, so one can get a
harmonic map by the well-known theorem of Eells-Sampson for the
existence of harmonic maps from this torus to $N$, which induces
an isomorphism from $\pi_1(T)$ to $\pi_1(N)$. Then a standard
argument shows that this harmonic map has constant energy density
by using the Bochner technique for harmonic maps. Furthermore, it
is a totally geodesic map and of real rank 1. Again, by means of a
result of Sampson \cite{s2}, this harmonic map maps $T$ to a
geodesic in $N$. This is a contradiction. This completes the proof
of the theorem 1.

\vskip 4mm \noindent {\bf Proof of Theorem 3:} Let $\Gamma$ be a
nonuniform lattice in $F_{4(-20)}$. (Without loss of generality,
one can assume that $\Gamma$ is torsion-free.) Assume it is the
fundamental group of a quasicompact K\"ahler manifold $M$. By
means of the Theorem 4, one has a harmonic map $f: M\to
H^2_{\mathbb{O}}\slash\Gamma$, which induces an isomorphism from
$\pi_1(M)$ to $\Gamma$. Here $H^2_{\mathbb{O}}$ is the Cayley
hyperbolic plane. Then a standard argument \cite{s1, y} shows that
$df(T^{1,0}_pM)$ ($p\in M$) can be regarded as an Abelian subspace
of the complexification ${\frak p}^{\mathbb{C}}$ of the tangent
space of $H^2_{\mathbb{O}}$. This tangent space can be identified
with the second factor of the Cartan decomposition $\frak g =
\frak t + \frak p$, here $\frak g$ is the Lie algebra of
$F_{4(-20)}$ and $\frak t$ is the Lie algebra of the maximal
compact subgroup in $F_{4(-20)}$. By the Lie-theoretic analysis in
\cite{ch}, one knows that the complex dimension of an Abelian
subspace in ${\frak p}^{\mathbb{C}}$ is at most $2$. So, we have
three cases to discuss for $df(T^{1,0}_pM)$.

{\bf (a):} ${\dim}_{\mathbb{C}}df(T^{1,0}_pM) = 1$, but
$df(T^{1,0}_pM)$ has real points.  So the real rank of $f$ is $1$.
By a result of J. H. Sampson \cite{s2}, $f$ maps $M$ to a geodesic
in $H^2_{\mathbb{O}}/\Gamma$. So, $\Gamma$ is isomorphic to the
ring of integers $\mathbb{Z}$. This is impossible;

{\bf (b):} ${\dim}_{\mathbb{C}}df(T^{1,0}_pM) = 1$, but
$df(T^{1,0}_pM)$ does not have real points. So the real rank of
$f$ is $2$. Similar to the previous proof, we have a decomposition
for $f$: $f = g\circ h$, here $h$ is a holomorphic map from $M$ to
an irreducible complex space $S$ of complex dimension $1$ and $g$
is a harmonic map. (Note that the discussion of Theorem 7.1 in
\cite{ct} does not work any more in the present case.) So, the
cohomological dimension of $\pi_1$, and hence $\Gamma$ is at most
$2$. But, the cohomological dimension of $\Gamma$ is actually $15$
(see \cite{b}). This is a contradiction;

{\bf (c):} ${\dim}_{\mathbb{C}}df(T^{1,0}_pM) = 2$. Completely
similar to the discussion in \cite{ch}, one also has a
decomposition for $f$: $f = g\circ h$, here $h$ is a holomorphic
map from $M$ to a quotient of the two-ball of finite volume and
$g$ is a geodesic immersion. Then, the same cohomological
dimension arguments show this is also impossible. This completes
the proof of the theorem.

\section{Examples of negatively $\delta$-pinched manifolds which
are not hyperbolic}

In this section, we will give some examples of negatively
$\delta$-pinched manifolds of finite volume which admit no
hyperbolic metric with finite volume under any smooth structure.
The constructing method is from Gromov and Thurston's paper
\cite{gt}. Actually, Gromov and Thurston constructed some examples
of compact manifolds which admit some $\delta$-pinched metric, but
no hyperbolic metric under any smooth structure. We only give a
sketch here. For detailed constructions, the reader can refer to
\cite{gt}.

Consider a non-singular quadratic form $\Phi_3$ in $4$ variables
$x_2, x_3, x_4, x_5$ with coefficients in ${\mathbb{Q}}$ and real
type $(1, 1, 1, -1)$. Let $\Gamma (\Phi_3)$ be the group of
automorphisms of the form $\Phi_3$ over the ring of integers
${\mathbb{Z}}$. It is well-known that $\Gamma (\Phi_3)$ may be
both cocompact and non-cocompact. We assume that $\Gamma(\Phi_3)$
is cocompact. Set the quadratic forms $\Phi_4 = (x_1)^2 + \Phi_3$
and $\Phi_5 = (x_0)^2 + \Phi_4$. Then one knows that
$\Gamma(\Phi_4)$ and $\Gamma(\Phi_5)$ are non-cocompact.
$\Gamma(\Phi_3)$ ($\Gamma(\Phi_4)$ respectively) can be considered
as a subgroup of $\Gamma(\Phi_4)$ ($\Gamma(\Phi_5)$ respectively)
and $H^3/(\Gamma(\Phi_3))$ is a compact totally geodesic
hypersurface of $H^4/(\Gamma(\Phi_4))$, while
$H^4/(\Gamma(\Phi_4))$ is a noncompact totally geodesic
hypersurface of $H^5/(\Gamma(\Phi_5))$. Then, Gromov-Thurston's
argument (Lemma 1.2 of \cite{gt}) can be applied to the space
pairs $(H^5/(\Gamma(\Phi_5)), D)$ and $(H^4/(\Gamma(\Phi_4)),
H^3/(\Gamma(\Phi_3)))$, here $D$ is a compact subset of
$H^4/(\Gamma(\Phi_4))$ containing $H^3/(\Gamma(\Phi_3))$. All
these together give the following
\begin{thm}
For every $\rho$, there are an orientable $5$-dimensional complete
manifold $V$ of constant curvature $-1$ and finite volume and a
compact totally geodesic orientable submanifold $V'$ of
codimension $2$ in $V$ such that the normal injectivity radius of
$V'$ in $V$ is greater than $\rho$; the corresponding homological
class of $V'$ in $V$ is trivial.
\end{thm}

Using the above constructed manifolds $V, V'$, one can get a
$Z_i$-ramified covering $\tilde{V}_i$ of $V$ at $V'$ for any
positive integer $i$. Then, on all these ramified coverings
$\tilde{V}_i$ with the induced smooth structure, one can construct
some complete negative curvature metrics $\tilde{g}_i$ whose
curvature at infinity is constant $-1$ and of finite volume since
$V'$ is compact. Now, consider all $\tilde{V}_i$ as topological
manifolds and suppose that there were complete metrics
$\tilde{g}^i$ of constant curvature $K=-1$ and finite volume on
each $\tilde{V}_i$ under some smooth structure of $\tilde{V}_i$
(not necessarily the above induced smooth structure). The action
of the cyclic group $Z_i$ on $\tilde{V}_i$ can be considered as a
deck transformation action. Because the dimension of $\tilde{V}_i$
is $5$, by means of the Mostow-Prasad rigidity theorem due to G.
Prasad \cite{p} (for nonuniform lattices in real hyperbolic
spaces), there exists an isometric action of $Z_i$ on
$\tilde{V}_i$ whose fixed point set $V''$ is homeomorphic to $V'$
and whose quotient $\overline{V}_i = \tilde{V}_i/Z_i$ has a
natural orbifold structure with constant curvature $K=-1$ and
finite volume. Also note that no two orbifolds $\overline{V}_i$
are isometric. On the other hand, by means of the volume estimate
argument in \cite{gt}, the volumes of $\overline{V}_i$ have a
uniform bound independent of $i$. So, Wang's finiteness theorem
for locally symmetric orbifolds (see \cite{w}, Theorem 8.1), which
asserts that there are at most finitely many isometric classes of
$n$$(\ge 4)$-dimensional complete orbifolds $\overline{V}$ with
$K(\overline{V}) = -1$ and $Vol(\overline{V}) \le$ a fixed
constant, implies a contradiction. Namely, there exists a positive
integer $i_0$, such that for $i \ge i_0$, $\tilde{V}_i$ does not
admit a complete metric of constant curvature $-1$ and finite
volume under any smooth structure. Finally, using the arguments of
$\S 3$ in \cite{gt}, one can show that there exist some
$\tilde{V}_i$ $(i\ge i_0)$ which admit no complete metric of
constant curvature $K = -1$ and finite volume under any smooth
structure, but carry some complete metrics with curvature $-1 \le
K \le -1-\epsilon$ and finite volume under the induced smooth
structure on $\tilde{V}_i$ from the ramified covering of $V$ at
$V'$ for arbitrary small positive $\epsilon$. Thus, combining
Gromov and Thurston's examples with ours, we actually obtain the
following
\begin{thm}
There exist some complete $5$-dimensional topological manifolds
which admit no complete metric of constant curvature $-1$ and
finite volume under any smooth structure, but carry some complete
metrics with curvature $-1\le K < -{\frac 1 4}$ and with finite
volume under some smooth structure. In other words, there exist
some groups which can be $\pi_1$ of some complete (open or closed)
$5$-dimensional negatively $\delta (>{\frac 1 4})$-pinched
manifolds of finite volume, but not a (uniform or nonuniform)
lattice of $SO(5, 1)$.
\end{thm}

\noindent {\bf Remarks:} A simple homotopical (or cohomological
dimension) argument shows that the groups in the above theorem
also cannot be (uniform or nonuniform) lattices of $SO(m, 1)$ for
any $m \neq 5$; in addition, these examples also show that Theorem
1 is a nontrivial generalization of the theorem in \cite{y}. A
natural problem is how to characterize these groups in algebraic
terms.

\bigskip
\noindent
J\"urgen Jost:\\
Max-Planck-Institute for Mathematics in the Sciences, Leipzig, Germany\\
and\\
\noindent
Yi-Hu Yang:\\
Department of Applied Mathematics, Tongji University, Shanghai, China\\
{\it e-mail}: yhyang@mail.tongji.edu.cn

\end{document}